\documentclass[12pt,a4paper]{article}

\usepackage[cp1251]{inputenc}
\usepackage[russian]{babel}%

\usepackage[centertags]{amsmath}
\usepackage{amsfonts}
\usepackage{amssymb}
\usepackage{graphicx}
\usepackage{amsthm}
\begin{document}




\oddsidemargin=-0.1cm
\topmargin=-5mm

\newtheorem{lemma}{Лемма}[section]
\newtheorem{remark}{Замечание}[section]
\newtheorem{theorem}{Теорема}[section]
\newtheorem{statement}{Section}[section]
\newtheorem{definition}{Определение}[section]
\newtheorem{proposition}{Предложение}[section]
\newtheorem{colorrary}{Следствие}[section]
\newtheorem*{theorem*}{Threorem}
\newtheorem*{theorem*1}{Theorem \normalfont{on polyhedra}}
\newtheorem*{corollary}{Следствие}
\newtheorem*{theorem*2}{Theorem on index}
\newtheorem*{conj}{Conjecture}
\def\proof{\par\noindent{\bf Доказательство}. \ignorespaces}
\def\endproof{{\ \vbox{\hrule\hbox{%
   \vrule height1.3ex\hskip1.3ex\vrule}\hrule
  }}\vspace{2mm}\par}

\def\proofff{\par\noindent{\bf Proof} \ignorespaces}
\def\endproofff{{\ \vbox{\hrule\hbox{%
   \vrule height1.3ex\hskip1.3ex\vrule}\hrule
  }}\vspace{2mm}\par}

\date{   }
\title{\bf Delone Sets:  Local Identity \\ and Global Symmetry
 \footnote{work is supported   by the RNF grant 14-11-00414}}
\author{Nikolay Dolbilin \\ Steklov Mathematical Institute \\ of the Russian Academy of Sciences\\
dolbilin@mi.ras.ru}


\maketitle

\renewcommand{\abstractname}{Abstract}
\renewcommand{\refname}{References}

 \hfill{\emph{To Friends who love and do
Geometry}}

\vskip 0.3cm

\begin{abstract}
In the paper we present a proof of the local criterion for crystalline structures which generalizes the
local criterion for regular systems. A Delone set is called a crystal if it is invariant with respect
to a crystallgraphic group.
So-called locally antipodal Delone sets,
i.e. such sets in which all $2R$-clusters are centrally symmetrical, are considered.
It turns out that the local antipodal sets have crystalline structure.
 Moreover, if in a locally antipodal set  all $2R$-clusters are the same
 the set is a regular system, i.e. a Delone set whose symmetry
 group operates transitively on the set.

\end{abstract}

\section{Introduction}

This  paper continues an investigative line started in the
pioneering work [5] on local conditions in a Delone set $X$  to
imply that   set $X$  is either a regular system, i.e. a
crystallographic orbit of a single point, or a crystal, i.e. the
orbit of a few points. On  Fig. 1 one can see the  set $X_1$ of
blue points which sit at nodes  of the square grid and the set
$X_2$ of red point quadruples. Each of these sets is a regular
system because each of them is an orbit of a 2D-crystallographic
group p4m (the full group of the standard square grid on the
plane). The union $X=X_1\cup X_2$ of the sets $X_1$ and $X_2$ is a
crystallographic orbit of two points, i.e. it is an  example of a
crystal.

So, a mathematical model of an ideal crystal uses two concepts:
the concept of the Delone set (which is of local character) and
the concept of the crystallographic group (which is of global
character). After Fedorov [8] a mathematical model of a
mono-crystalline matter is defined  as a Delone set which is
invariant with respect to some crystallographic group. One should
emphasize that under this definition the well-known periodicity of
crystal in all 3 dimensions is not an additional requirement. By
the famous Sch\"onflies-Bieberbach theorem [3,4], any space group
contains a translational subgroup with a finite index.

Since the crystallization is such a process that results from
mutual interaction of  nearby atoms, it is believed (R.~Feynman,
N.V.~Belov, at al) that the long-range order of atomic structures
of crystals (and quasi-crystals too) comes out of local rules
restricting the arrangement of nearby atoms. R.~Feynman wrote
([10], Ch. 30): "when the atoms of matter are not moving around
very much, they get stuck together and arrange themselves in a
configuration with as low an energy as possible. If the atoms in a
certain place have found a pattern which seems to be of low
energy, then the atoms somewhere else will probably make the same
arrangement. For these reasons, we have in a solid material a
repetitive pattern of atoms. In other words, the conditions in a
crystal are this way: The environment of a particular atom in a
crystal has a certain arrangement, and if you look at the same
kind of an atom at another place farther along, you will find one
whose surroundings are exactly the same. If you pick an atom
farther along by the same distance, you will find the conditions
exactly the same once more. The pattern is repeated over and over
again, and, of course, in three dimensions." The crystallographer
N.V.~Belov also suggested similar arguments in a problem "on the
501st element".

However, before 1970's there were no whatever rigorous results to explain a
link between properties of local patterns
and the global order in the internal  structure of crystals    until
Delone and his students initiated developing the  local theory of crystals [5].
One of two  main aims of the local  theory was (and is) rigorous derivation of space
group symmetry of a crystalline structure from the pairwise identity
of local arrangements around each atom.

One should mention that the link between the identity of local fragments of a structure
and global order of the structure
seemed  obvious,  and searching for  an exact wording of  this connection and its rigorous proof were

 seen  a purely abstract goal
that would be of interest only to mathematicians.

However, the subsequent discovery of non-periodic Penrose patterns
(1977)  and the discovery by D.~Shechtman  of real quasicrystals
(1982, Nobel Prize in 2011) have showed that there are
non-periodic structures which are as locally identical  as
crystals. These discoveries suggest that the connection  between
the local identity and the global order is not so obvious. One of
the goals of the local theory of global order was
 to look for right wordings of statements  and then prove them.

The local theory has been developed for Delone sets as well as for
polyhedral tilings, including combinatorial aspects of the theory
(see e.g. [11,12,13]).

The paper is organized as follows. In the next section we give
definitions of all necessary concepts  and  short
 survey of some earlier results. Then we give formulation of  the local criterion for a crystal and of several  new
 "local"
 theorems on locally antipodal Delone sets (Theorems 1-- 5),  which will be proved in
concluding sections of the paper.

\section{Basic Definitions  and Results}

\textbf{Definition 2.1}.   A point set  $X\subset \mathbb{R}^d$ is
called a \emph{Delone set} with parameters  $r$ and  $R$, where
$r$ and $R>0$, (or an  $(r,R)$-\emph{system}, see [1,2] ), if two
conditions hold:
\newline (1) an open  $d$-ball  $B^o_y(r)$ of radius $r$ centered at an  \emph{arbitrary} point  $y\in \mathbb{R}^d$
   contains at most  one point from  $X$:
$$|B^o_y(r)\cap X|\leq 1; \eqno(r)$$
(2) a closed  $d$-ball  $B_y(R)$ of radius $R$ centered at an
arbitrary point $y$  contains at least one point from  $X$:
$$|B_y(R)\cap X|\geq 1. \eqno(R)$$
We note that by condition  $(r)$ the distance between any two
points $x$ and $x'\in X$ is not less than $2r$.

For $x\in X$ we denote $C_x(\rho):=X\cap B_x(\rho)$  and will call
the set  $C_x(\rho)\subset X$  a  $\rho$-\emph{cluster} of point
$x$. Thus, a $\rho$-cluster $C_x(\rho)$ consists of all points of
$X$  that are placed from $x$ at distance at most $\rho$. It is
easy to see that for  $\rho<2r$ $C_x(\rho)=\{x\}$. It is
well-known that for   $\rho\geq 2R$, the $\rho$-cluster
$C_x(\rho)$ of any point $x\in X$ has the full rank: the dimension
of conv($C_x(\rho))=d$

In principle, the $\rho$-cluster $C_x(\rho)$ is considered as a
pair (the center $x$,  the point set $C_x(\rho)$). However, since
notation $C_x(\rho)$ contains information on the center $x$ we can
miss the notation of a cluster as the pair.  We emphasize that we
distinguish between
 $\rho$-clusters  $C_x(\rho)$ and  $C_{x'}(\rho)$ of different points
 $x$ and
$x'$, even if the two sets  generally may coincide (see Fig. 1).

\begin{center}
\begin{figure}[!ht]
\hbox{\includegraphics[width=8cm]{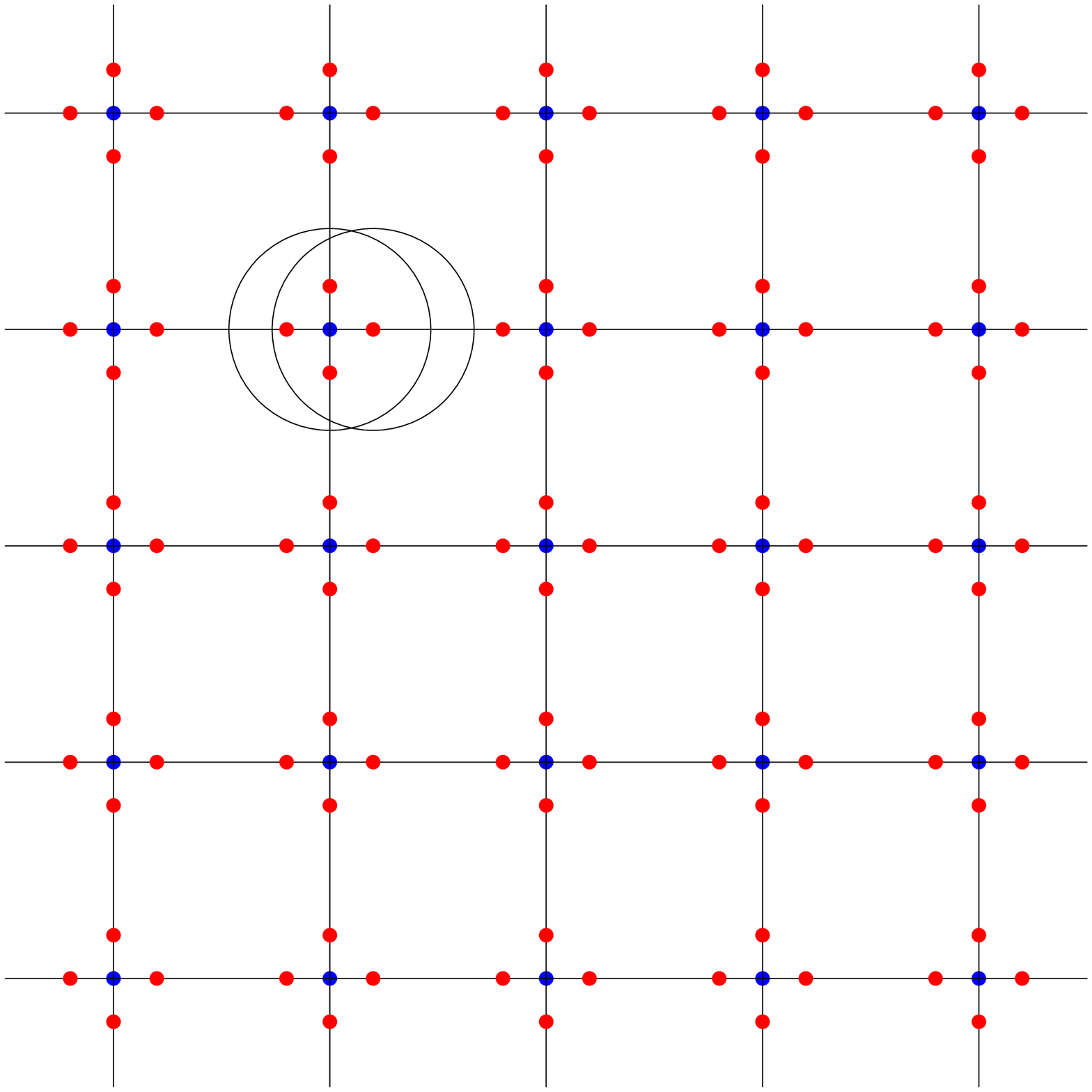}
\raise5.6cm\llap{$x$\kern5.1cm}
\raise5.9cm\llap{$x'$\kern4.6cm}
\raise0.0cm\llap{Fig. 1\kern3.5cm} }
\end{figure}
\end{center}

\textbf{Definition 2.2}.   Two  $\rho$-clusters  $C_x(\rho)$ and
$C_{x'}(\rho)$ are called \emph{equivalent}, if there is an isometry
$g\in O(d)$ such that
\newline\centerline{$g:x\mapsto x'$ and}
\newline\centerline{ $g: C_x(\rho)\rightarrow C_{x'}(\rho)$.}

We emphasize that the requirement of equivalence of two clusters
is some stronger than  just of congruence of sets of points that
enter these clusters. Two clusters depicted on fig.1 around two
points $x$ and $x'$ coincide as  subsets of $X$. However, since
this subset of $X$ surrounds the points $x$ and $x'$ in different
ways it is natural to distinguish between the  $\rho$-clusters
$C_x(\rho)$ and $C_{x'}(\rho)$.  Indeed, the clusters $C_x(\rho)$
and $C_{x'}(\rho)$ are non-equivalent because there is no isometry
that moves both point $x$ and cluster $C_x(\rho)$ into point $x'$
and cluster $C_{x'}(\rho)$, respectively.

In a Delone set for any $\rho>0$ a set of clusters is partitioned
into classes of equivalent $\rho$-clusters. For any given $\rho$
if $\rho<2r$,  $\rho$-cluster at any point of $X$ consists of a
single point: $C_x(\rho)=\{x\}$, i.e. all "small" $\rho$-clusters
in $X$ are equivalent. Given Delone set $X$, we denote by
$N(\rho)$  the cardinality of a set of equivalence classes of
$\rho$-clusters in  $X$.

We have in any Delone set $X$ $N(\rho)=1$ for $\rho<2r$. However,
for larger $\rho$, $\rho > 2r$, $N(\rho)$ can be generally
infinite.

\textbf{Definition 2.3}.
 A Delone set
$X$ is said to be \emph{of finite type} if for each  $\rho>0$ the
number  $N(\rho)$ of classes of $\rho$-clusters is finite.

As we said, for any Delone set $X$ function $N(\rho)$ is always
defined and equal to 1 for all $\rho<1$. It is not hard to prove
the following:

\smallskip
\textbf{Statement 2.1}. \emph{Function $N(\rho)<\infty$ for all $\rho>0$ if
$N(2R)< \infty $.}

\smallskip
The key reason of this fact  is as follows. Given $X$, the
condition $N(2R)<\infty$ implies, that in a Delone tiling for set
$X$ (see [2]) there are just finitely many pairwise non-congruent
tiles. Next, we note that in the Delone  tiling, that is,
importantly, face-to-face, any two convex finite $d$-polyhedra $P$
and $Q$ that share a common $(d-1)$-face can form  just finitely
many non-congruent pairs $(P,Q)$. From here it follows in the
Delone tiling with $N(2R)<\infty $ for any $\rho$ there are just
finitely many non-congruent fillings of a ball $B_x(\rho)$. The
finiteness of different parts of the tiling of size $\rho$ implies
the finiteness of $N(\rho)$.

Now we take a point $x\in X$ and the cluster  $C_x(2R)$. Points of
the cluster uniquely determine all Delone cells  for $X$ that meet
at point $x$. Now, since $N(2R)<\infty$ the number of
non-congruent Delone cells in a Delone tiling for $X$ is finite.
Due to this finiteness and also to above-mentioned  finiteness of
the number of pairs of $P$ and $Q$ glue along their common
hyperface. These two sorts of finiteness imply just
 finite number of non-congruent  fillings of a ball of radius
 $\rho$ with Delone tiles. It follows that every
 $2R$-cluster $C_x(2R)$ admits just a finite number of non-congruent extensions
 to a $\rho$-cluster $C_x(\rho)$. Since, by assumption,  there are finitely many pairwise non-congruent
  $2R$-clusters,  then there are finitely many  $\rho$-clusters  for any given  $\rho>0$.

\medskip
From now on, we will consider Delone sets of finite type only. We
note that in such a Delone set the number  $N(\rho)$ of
$\rho$-clusters is a positive, integer-valued, non-decreasing,
piece-wise constant  function, continuous on the right.

Very important examples of Delone sets are the so-called
\emph{regular systems} and \emph{crystals}. The concept of the
regular system was studied by E.S.Fedorov [8]. Regular systems are
 discrete homogenous sets which
 looks the same up to infinity from any its point.  Here is an equivalent
definition in terms of a Delone set.

\medskip

\textbf{Definition 2.4}  A \emph{regular system}  is  a Delone set $X\subset \mathbb{R}^d$
whose symmetry group acts transitively, i.e. for any two points
$x$ and $x'\in X $ there is an isometry  $g\in Iso(d)$ such that
\newline\centerline{$g:x\mapsto x'$ и $g: X\rightarrow X$}.

\medskip

Recall that a group  $G\subset Iso(d)$ is called a
\emph{crystallographic group}  if
\newline (1)  $G$ operates discontinuously at each point  $y\in \mathbb{R}^d$, i.e.
if for any point $y\in \mathbb{R}^d$ orbit $G\cdot y$ is a
discrete set;
\newline (2) $G$ has a compact fundamental domain.

One can reformulate the definition of a regular set in terms of a
crystallographic group. We emphasize that in the following
statement we do not require that $X$ is a Delone set. This
condition results from properties of a crystallographic group.

\medskip
\noindent\textbf{Statement 2.2}. \emph{ A set $X\subset
\mathbb{R}^d$ is a regular system if and only if the set $X$ is an
orbit of a point  $x\in \mathbb{R}^d$ with respect to a
crystallographic group  $G\subset Iso(d)$.}

\smallskip

A regular set is an important particular case of the more general
concept of a \emph{crystal}.

\smallskip
\textbf{Definition 2.5}.  A \emph{crystal} is a Delone set $X$
such that $X$ is a finite collection of orbits with respect to its
symmetry group Sym$(X)$: $X=$Sym$(X)\cdot X_0$, where $X_0$ is a
finite point set.

It is not hard to prove  that the  symmetry group of a crystal is
a crystallographic group. Thus we have the following  statement.

\smallskip
 \textbf{Statement 2.3}. \emph{
 A set $X\subset \mathbb{R}^d$ is a crystal if and only if it is an orbit of
 a finite set $X_0$ with respect to a crystallographic group
 $G$, i.e.  $X=G\cdot X_0$.}

\smallskip
 These classes of Delone sets can be described via the cluster counting function  $N(\rho)$ as follows. A Delone set
 of finite type  is a regular system if and only if
$N(\rho)\equiv 1$ on  $R_+$. A Delone set is a crystal if and only
if its cluster counting function is bounded:
$$
N(\rho)\leq m <\infty, \,\,\, \hbox{ where } m\leq |X_0|.
$$
If  $m=1$, then a crystal is a regular system.

\smallskip
 The above-mentioned  definitions of a regular system and a crystal go back to
 Fedorov, [8]. Earlier, before Fedorov's work, a crystal had been  considered as a set of pairwise congruent and parallel
 lattices. The definition of a crystal in terms of regular systems
 seemed to generalize the Haui-Bravais concept of a crystal as a 3D periodic Delone set.

 However, indeed, due to the  Sch\"onflis-Bieberbach theorem, the  more general structure of a regular system
 is also  the union of lattices. Therefore, due to the
 Fedorov definition, a crystal is  the union of several lattices exactly as in the Haui-Bravais approach.

\bigskip
Indeed, let $h$ be the index of the translational subgroup of a
crystallographic group $G\subset Iso(d)$, $|X_0|$ the number of
points in $X_0$, then  a crystal $G\cdot X_0$ splits into $m$
pairwise congruent and parallel lattices of rank $d$:
$$
G\cdot X_0= \cup_i^m (T\cdot x_i \cup T \cdot g_2(x_i) \cup \ldots
\cup T\cdot g_h(x_i)), \,\,\, x_i\in X_0.
$$
We note that $m$ is strictly smaller than $h\cdot |X_0|$ if some
 points $x_i$ and $x_j$ in $X_0$ belong to one $G$-orbit.

Now we define the \emph{group} $S_x(\rho)$ of $\rho$-cluster
$C_x(\rho)$ as a subgroup of Iso$(d)$ to consist of those
isometries $s$, such that
$$s:x\mapsto x, \,  s: C_x(\rho)\longrightarrow C_x(\rho).$$
\noindent Let us denote by $M_x(\rho)$ the order of the group
$S_x(\rho)$. Since the rank of $C_x(2R)$  in a Delone set
$X\subset \mathbb{R}^d$  equals $d$, the order $1\leq M_x(\rho)
<\infty$, i.e. the  function $M_x(\rho)$ is defined for all
$\rho\geq 2R$.

The function $M_x(\rho)$ for all $\rho\geq 2R$ takes positive
integer values, is continuous on the left and  non-increasing.
Moreover,  the ratio $M_x(\rho):M_x(\rho')$ is integer if
$\rho'>\rho$. In fact, it is obvious that the group $S_x(\rho')$
of a bigger cluster
 $C_x(\rho')$  either coincides  with  $S_x(\rho)$, or it is a proper subgroup
of $S_x(\rho)$,  $\rho'>\rho$.

Let  $X$ be a Delone set of finite type. Then for any positive $\rho $ the  set $X$ splits into a finite number $N(\rho)$ of subsets
 $X_1, X_2, \ldots , X_{N(\rho)}$, such that points $x$ and  $x'$ from every one subset  $X_i$ have equivalent
  $\rho$-clusters $C_x(\rho)$ и
$C_{x'}(\rho)$, but at points from different subsets  $X_i$ and
$X_j$ the $\rho$-clusters are not equivalent. The groups of
equivalent  $\rho$-clusters are conjugate in  $Iso(d)$ and
consequently have the same order  $M_i(\rho)$, where  $i$ is the
index of a subset $X_i$, $i\in
 [1,N(\rho)]$.

\bigskip
One of main goals of the local theory of regular systems is to
determine a radius $\hat \rho$ such that any Delone set $X$ (with
parameters $r$ and $R$) with $N(\hat \rho)=1$ is a regular system.
Certainly, the answer may depend on the dimension. So, for $d=1$
it is easy to see that a Delone set on a line is a regular system
if $N(2R)=1$. The value $2R$ cannot be improved: in fact, for any
$\varepsilon>0$ there are Delone sets with $N(2R-\varepsilon)=1$
that are not regular systems. The first important  result in  the
local theory of regular systems was obtained in [5].

\bigskip
 \textbf{Theorem 2.1} (Local criterion for regular systems).
 \emph{A Delone set  $X\subset \mathbb{R}^d$ (with parameters  $r$ and  $R$
is a regular system if and only if for some  $\rho_0>0$ the
following  conditions hold:}
\newline (I) $N(\rho_0+2R)=1$;
\newline (II) $M(\rho_0)=M(\rho_0+2R)$.

\medskip

Condition (I) means that   $(\rho_0+2R)$-clusters at all points
$x\in X$ are equivalent. Therefore the groups $S_x(\rho_o+2R)$ of
the  clusters are pairwise conjugate. Due to condition (II) for
each point $x\in X$ the groups  $S_x(\rho_0)$ and $S_x(\rho_0+2R)$
coincide.

\medskip Let us select among Delone sets with $N(2R)=1$
\emph{locally asymmetric}  sets, i.e. such that the group
$S_x(2R)$ is trivial. From the local criterion follows a theorem
for locally asymmetric Delone sets.

\bigskip
\textbf{Theorem 2.2} [Locally asymmetric sets]. \emph{Let a Delone set
$X\subset\mathbb{R}^d$ be locally asymmetric set and $N(4R)=1$.
Then $X$ is a regular system, i.e. $N(\rho)\equiv 1, \forall
\rho>2R$.}

\medskip

The theorem immediately follows from the criterion for
$\rho_0=2R$.

It is amazing that due to the following theorem
 the condition $N(4R)=1$ can not be reduced.

\bigskip
 \textbf{Theorem 2.3} ($(4R-\varepsilon)$-theorem ). \emph{For any given  $\varepsilon>0$
and any dimension $d$, $d\geq 2$ there is a Delone set  $X\subset
\mathbb{R}^d$ such that  $N(4R-\varepsilon)=1$, but $X$ is not a
regular system.}

\medskip The theorem is proved by means of an explicit
construction.
Below we present such a design for dimension $d=2$. This
construction can be  easily extended to any dimension $d$.

We start  constructing the design  with a rectangular lattice
$\Lambda$ (Fig. 2, on the left). The lattice $\Lambda$ has a
fundamental rectangle with side length $a$ and $b$ where $a<<b$.
It is clear that the parameter $R=\sqrt{b^2+a^2}/2$. Since $a<<b$
we have  $$2R\sim b(1+a/2b)=b+a/2.$$

Horizontal rows of the $\Lambda$ form a bi-infinite sequence with
indices $i\in \mathbb{Z}$. The set of the rows splits into couples
$P_{2i}$ of rows with sequel indices $(i,i+1)$ where the first one
$i$ is even. Let us take $c$ so that $0<c<a/2$ and shift
\textbf{each}  couple $P_{i+1}$ relatively to the previous couple
$P_i$ by $c$ to the left or to the right.

We get a sequence of mutually shifted rows that can be encoded a
bi-infinite sequence $l=\ldots RLLRL\ldots$. There are uncountably
many different bi-infinite binary sequences $\{l\}$.    The
corresponding Delone sets $X_l$ have  the same parameters $r$ and
$R$.  Among the sequences $\{l\}$ there are exactly 3 ones such
that the corresponding Delone sets are regular systems.  Two
sequences $\ldots LLLL \ldots $ and  $\ldots RRRR \ldots $
generate two congruent each other  regular systems. The third
bi-infinite sequence $\ldots RLRLRL\ldots $ encodes one more
regular system to be mirror symmetrical to itself. All the other
Delone sets from the family are not regular systems, though as it
is easy to see that all of them have the same $b$-clusters
$C_x(b)$. Since $b\sim 2R-a/2$ and $a>0$ can be chosen arbitrarily
small, we get the theorem.

\begin{center}
\begin{figure}
\null \hspace{4.6cm} $\cdots RRR\cdots$ \hspace{1.0cm} $\cdots
RLRL\cdots$ \hspace{0.5cm} $\cdots RLL\cdots $ \vspace{1.0cm}

\resizebox{3.4cm}{!}{\includegraphics{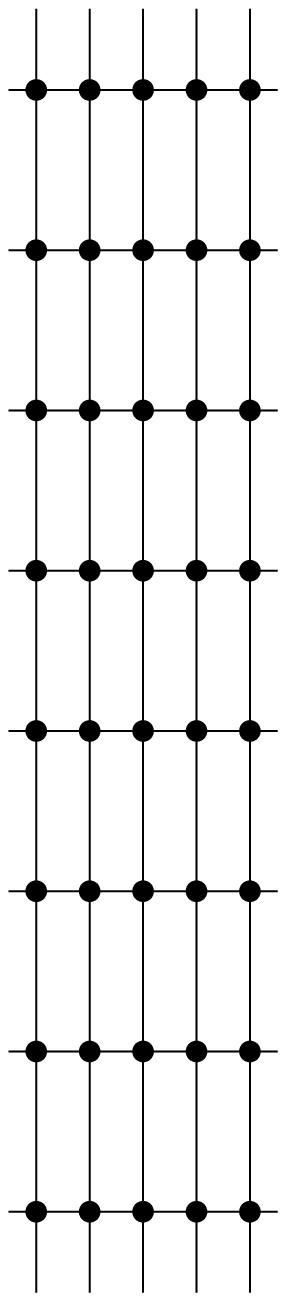}} \hspace{-1.5cm}
 \resizebox{4.0cm}{!}{\includegraphics{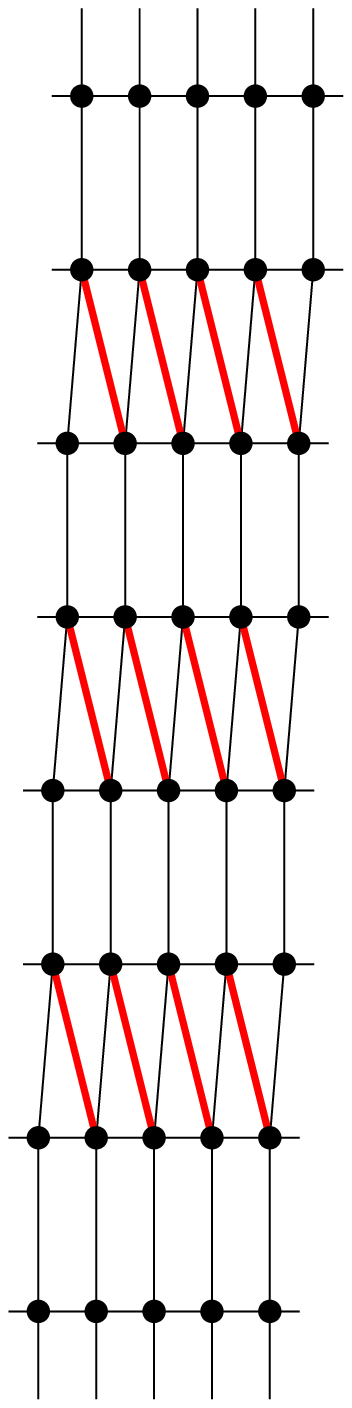}} \hspace{-2cm}
\resizebox{4.0cm}{!}{\includegraphics{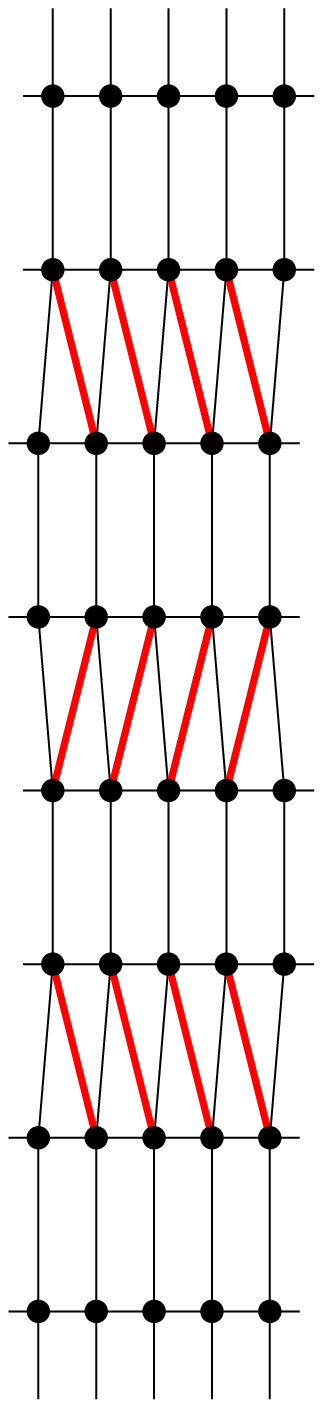}} \hspace{-1cm}%
\resizebox{4.0cm}{!}{\includegraphics{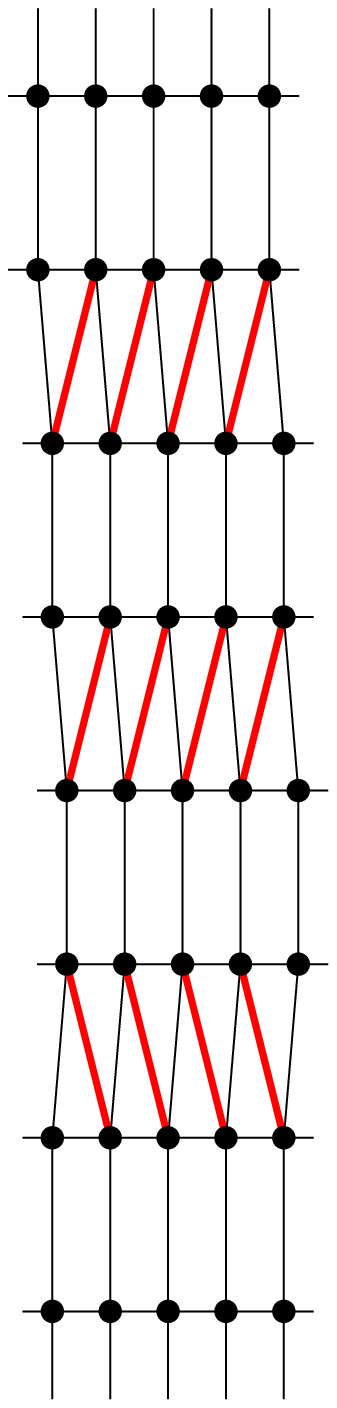}}
\centerline{Fig.2}
\end{figure}
\end{center}


In this context it is particularly interesting to note that there
are the following results  for Delone sets in an Euclidean plane
and in 3D space:

\medskip
\textbf{Theorem 2.5} (Regular systems, $d=2,3$).

\noindent  (1) \emph{Let  $X\subset \mathbb{R}^2$ be a Delone set
in plane, if $N(4R)=1,$ then $X$ is a regular system.}

\noindent (2) \emph{Let and $X\subset \mathbb{R}^3$ be a Delone
set, if $N(10R)=1,$ then $X$is a regular system. }

\medskip This result was obtained by M.~Stogrin and by N.~
Dolbilin independently  years ago; however  a (complete) proof
remains unpublished if one does not take into account publications
on key ideas. As for  point (1) of theorem 2.5, case $d=2$, here
we just mention that the part of  the theorem can be derived from
 the following  theorem:

\medskip
\textbf{Theorem 2.6} [15]. \emph{A tiling of Euclidean plane by
convex polygons is regular, i.e. a tiling with a transitive
symmetry group, if all first coronas are equivalent}.

\medskip

Details of the proof  of the $10R$-theorem for regular systems in
$\mathbb{R}^3$ will appear in [14]. Here we just mention that, due
to the $(4R-\varepsilon)$-theorem,  the estimate $4R$ for plane is
the best estimate. As for the estimate $10R$ for 3D-space, it
looks much higher than the actual one. The difficulty lies in the
fact that we can not deal effectively with the $4R$-cluster group.

\bigskip

In this regard, it is especially remarkable that in the case where
 the $2R$-cluster group contains the central symmetry, a
sufficient condition on regular systems becomes extremely simple
and holds for any dimension.

\textbf{Definition 2.6}. A Delone set $X$ is called a
\emph{locally antipodal} set if
 the $2R$-cluster $C_x(2R)$ for each point $x\in X$ is centrally symmetrical
 about the cluster center $x$.

\medskip Now we present several theorems to be proved in the next
sections.

\medskip

\textbf{Theorem 1} . \emph{If $X$ is  a locally antipodal set and
$N(2R)=1$, then $X$ is a regular system.}

\medskip

\textbf{Theorem 2}. \emph{A locally antipodal Delone set $X$ is
centrally symmetrical about each point $x\in X$ globally.}

We emphasize that  in either  theorem 2 or in the following
theorem 3, no condition on the cluster counting function $N(\rho)$
is required. Moreover, we do not require even that a Delone set
$X$ is of finite type.

\textbf{Theorem 3}. \emph{A locally antipodal Delone set $X\subset
\mathbb{R}^d$ is a crystal. Moreover, $X$ is the disjoint union of
at most $2^d-1$ congruent and parallel  lattices:}
$$
X = \bigsqcup\limits_{i = 1}^n (x + \lambda _i/2 + \Lambda),
$$
where $\Lambda$ is the  maximum  lattice for $X$ such that
$X+\Lambda=X$, $\lambda_i\in \Lambda$ and \newline
$\lambda_i\equiv \lambda_j (\mod 2)\Leftrightarrow i=j$, $n<2^d$.

\medskip
Theorems 1-- 3 have been published in part in [7] and [8]. In this
paper theorems 1 and 2 are easily derived from  the following
theorem.

\medskip
\textbf{Theorem 4} (Uniqueness theorem). \emph{Let $X$ and $Y$ be
 Delone sets with the  parameter $R$,  suppose  they
 have a point $x$ in common. Let the  $2R$-clusters of $X$ and $Y$ centered at this point $x$ coincide, i.e.
   $C_x(2R)=C'_x(2R)$,  where $C_x(\rho)$ stands for a cluster in $X$ and
 $C'_y(\rho )$ for a cluster in $Y$. Then $X=Y$.}

\bigskip

In the conclusion to this  section we present  a local criterion
for a crystal. This criterion generalizes the local criterion for
regular systems. It was announced [6] and proved a while ago but a
full proof was published recently [7] (in Russian). The proof in
this paper is a slight improvement  of the proof in [7].

\textbf{Theorem 5} (Local criterion for a crystal). \emph{A Delone
set $X$ of finite type is a crystal which consists of $m$ regular
systems if and only if there is some $\rho_0>0$ such that two
conditions hold}:
\newline
1) $N(\rho_0)=N(\rho_0+2R)=m;$
\newline
 2) $S_x(\rho_0)=S_x(\rho_0+2R), \forall x\in X.$ .

\medskip
The local criterion for regular systems (theorem 2.1) is a
particular case of theorem 5. Indeed, the  condition
$N(\rho_0+2R)=1$ implies $N(\rho_0)=N(\rho_0+2R)=1$. The cluster
groups $S_x(\rho)$ for all $x\in X$ are pairwise conjugate, hence
it suffices to require $S_x(\rho_0)=S_x(\rho_0+2R)$ for one point
from $X$ only.


\section{Proof of the Local Criterion for a Crystal }

First of all we make some comments on conditions 1) and 2) of
theorem 5. Condition 1) means that with increasing radius $\rho$
the number of cluster classes on segment $[\rho_0,\rho_0+2R]$ does
not increase, i.e. remains unchanged: $N(\rho)= N(\rho+2R)$.

In addition, due to   condition  2), the   cluster group
$S_x(\rho_0)$, $\forall x\in X$, does not get smaller under the
$2R$-extension of $\rho_0$-cluster: $S_x(\rho_0)=S_x(\rho_0+2R)$.
The key point  of theorem 5 is that the stabilization of these two
parameters (the number of cluster classes and the order of cluster
groups) on segment $[\rho_0,\rho_0+2R]$ implies their
stabilization on the rest of the half-line  $[\rho_0+2R, \infty)$.

\begin{lemma}[on the $2R$-chain]
For any pair of points $x$ and  $x'\in X$, where  $X$ is a Delone
set, in $X$  there is a finite sequence  $x_1=x, x_2, \ldots,
x_k=x'$, such that  $|x_ix_{i+1}|<2R$ для $i\in [1, k-1].$
\end{lemma}

We omit a proof of the lemma which can be found in, e.g.   [5].

We recall that $X$ splits into disjoint subsets
$X=\bigsqcup_{i=1}^m X_i$, where $X_i$ is a subset of all points of $X$ whose $\rho_0$-clusters
are pairwise equivalent and belong to the $i$-th class.

\medskip
\begin{lemma}[on the  $2R$-extension]  Let  a Delone set $X$ fulfil conditions 1) and 2) of theorem $5$ and  $x$,  $x'\in
X_i$. Let  $f\in Iso(d)$ be  an isometry such that
$$
f:x\mapsto x' \hbox{ and } f: C_x(\rho_0)\rightarrow C_{x'}(\rho_0).
\eqno(1)
$$
Then the same isometry $f$ superposes the concentrically  bigger
cluster $C_x(\rho_0+2R)$  onto  cluster $C_{x'}(\rho_0+2R)$:
$$
f: C_x(\rho_0+2R)\rightarrow C_{x'}(\rho_0+2R). \eqno(2)$$
\end{lemma}

\textbf{Proof}.  If the $\rho$-clusters $C_x(\rho)$ and $C_{x'}(\rho)$
are equivalent, by condition 1)  of theorem 5, the corresponding
$(\rho_0+2R)$-clusters  are equivalent too. Therefore there is
 an isometry $g$ such that
$$
g: C_x(\rho_0+2R)\rightarrow C_{x'}(\rho_0+2R). $$
If $g=f$ the lemma is already proved. Now we assume $f\neq g $ we consider
the superposition of isometries  $f^{-1}\circ g$. The order here is
from the right to the left:  $g$ is applied first, followed by $f^{-1}$.
$$
f^{-1}(g(C_x(\rho_0)))= f^{-1}(C_{x'}(\rho_0))=C_x(\rho_0).
$$
We have:
$$f^{-1}\circ g:x\mapsto x \hbox{ and } f^{-1}\circ g: C_x(\rho_0)\rightarrow
C_x(\rho_0). \eqno(3)$$

Relationship  (3) shows that $f^{-1}\circ g$ is a symmetry $s$ of
the  $\rho$-cluster $C_x(\rho_0)$: $s\in S_x(\rho_0)$. By
condition 2) of theorem 5 we also have   $s\in S_x(\rho_0+2R)$.
The relation $f^{-1}\circ g=s$ implies $g\circ s^{-1}=f$.
Therefore one gets
$$f(C_x(\rho_0+2R))=(g\circ
s^{-1})(C_x(\rho_0+2R))=$$
$$= g(s^{-1}(C_x(\rho_0+2R)))=
g(C_x(\rho_0+2R))=C_{x'}(\rho_0+2R).$$\endproof

\begin{lemma} [Main Lemma].  Let a Delone set  $X$ fulfil conditions $1)$ and  $2)$ of theorem 5 and $X_i$ a subset of $X$
of  all points whose $\rho_0$-clusters belong to the $i$-th class,
$i\in [1,m]$. Let  a group $G_i=<f>$ be generated by all
isometries $f$ such that $f:C_x(\rho_0)\rightarrow
C_{x'}(\rho_0)$, $\forall x,x'\in X_i$. Then $G_i$
 operates transitively on every set $X_j$,
 $\forall j\in [1,m].$
Moreover, the group $G_i$ does not depend on $i$, $G_i=$Sym$(X)$.
\end{lemma}

\textbf{Proof}. In fact, since for $x$ and  $x'\in X_i$ the
$\rho_0$-clusters $C_x(\rho_0)$ and $C_{x'}(\rho_0)$ are
equivalent, there are  several  isometries superposing these
clusters.  The number of these isometries is equal to the order
$|S_x(\rho_0)|$ of the cluster group.

We prove that if $f$ is
 one of these  isometries, then $f$ is a symmetry of the whole $X$. First we take an
arbitrary point $y\in X$ and prove that its image $f(y)$ belongs
to $X$ (see Fig.3). Let us connect points  $x$ and  $y$ with a
$2R$-chain  $\cal L$:
$${\cal L}=\{x_1=x,x_2, \ldots , x_n=y \,: \,|x_ix_{i+1}|<2R, \forall i\in[1,n-1]|\}.
$$

\begin{center}
\begin{figure}[!ht]
\hbox{\includegraphics[width=8cm]{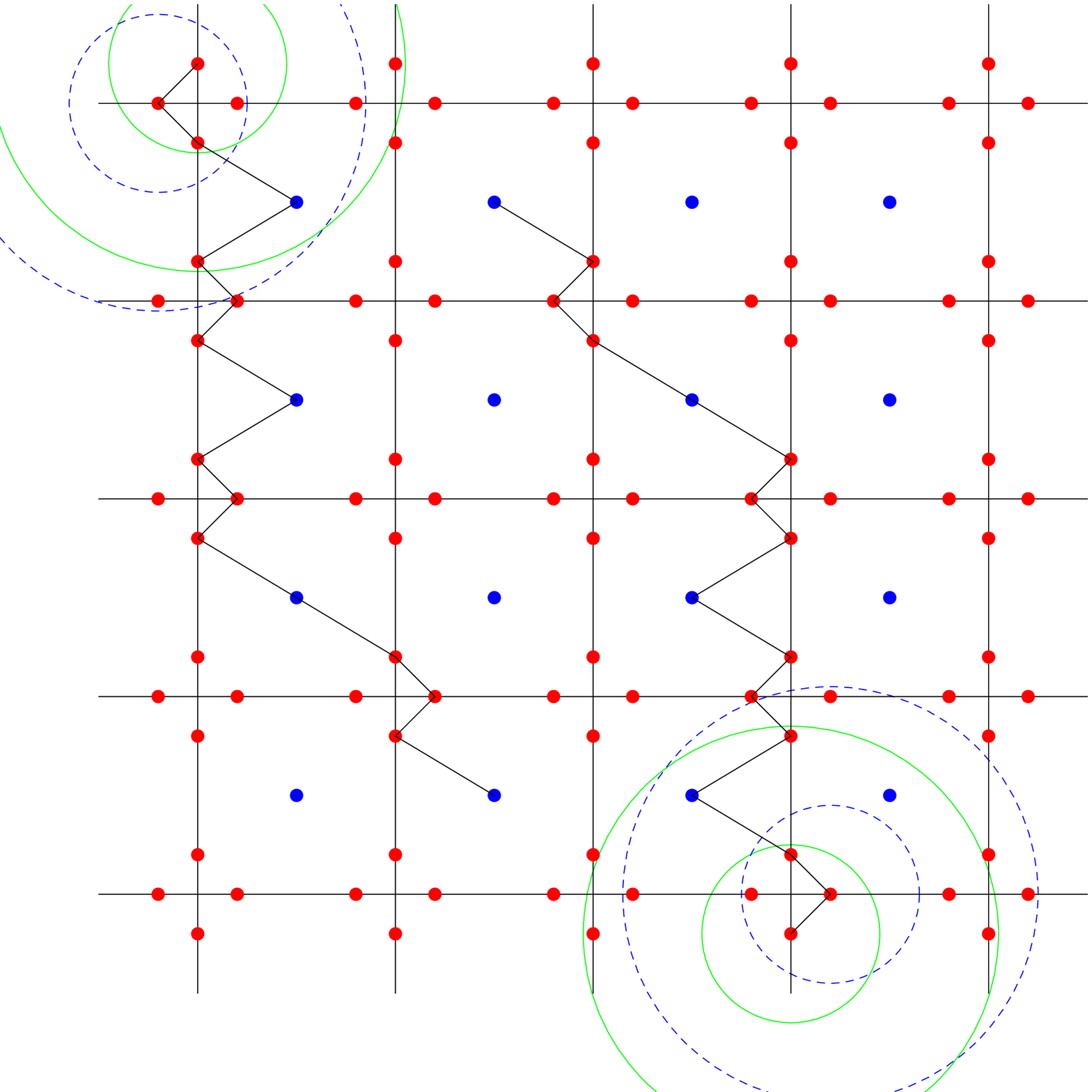}
\raise7.5cm\llap{$x$\kern6.4cm}%
\raise7.0cm\llap{$x_2$\kern6.9cm}%
\raise2.0cm\llap{$y$\kern4.6cm}%
\raise6.1cm\llap{$y'$\kern4.5cm}%
\raise1.0cm\llap{$x'$\kern2.0cm}%
\raise1.6cm\llap{$x_2'$\kern1.7cm}%
\raise2.3cm\llap{$x_{n-1}$\kern4.9cm}%
\raise6.0cm\llap{$x'_{n-1}$\kern2.9cm}%
\raise0.0cm\llap{Fig. 3\kern3.6cm} }
\end{figure}
\end{center}

Since $f(C_{x_1}(\rho_0))=C_{x'_1}(\rho_0)$, then by Lemma 3.2
$$f(C_{x_1}(\rho_0+2R))=C_{x'_1}(\rho_0+2R). \eqno(4) $$
Since  $|x_1x_2|<2R$ we have
$$C_{x_2}(\rho_0)\subset C_{x_1}(\rho_0+2R).$$
Therefore relation (4) implies:
$$f:C_{x_2}(\rho_0)\rightarrow C_{x'_2}(\rho_0).$$
By lemma  3.2 we have
$$f: C_{x_2}(\rho_0+2R)\rightarrow C_{x'_2}(\rho_0+2R). \eqno(5)$$
From  the inequality $|x_2x_3|<2R$ it follows that
$$C_{x_3}(\rho_0)\subset C_{x_2}(\rho_0+2R).$$
Therefore due to  (5) we have:
$$f:C_{x_3}(\rho_0)\rightarrow C_{x'_3}(\rho_0).$$
By Lemma 3.2 we have again:
$$f:C_{x_3}(\rho_0+2R)\rightarrow C_{x'_3}(\rho_0+2R).$$
Moving along chain  $\cal L$ and repeating this argument finitely
many times we get that the  $2R$-chain   ${\cal L} \subset X$ is
moved by  isometry $f$ into a  $2R$-chain  ${\cal L}'\subset X$.
The endpoint  $y$ of the first chain $\cal L$ moves into the
endpoint $y'$ of the second one. Thus, it is shown that the
isometry $f$ maps  $X$ \textbf{into} $X$: $f(X)\subseteq X$.

We show now that this map is a map \textbf{onto} the whole $X$:
$f(X)\supseteq X$. Let us take an arbitrary point  $y''\in X$ and
show that its pre-image $f^{-1}(y'')$ also belongs to $X$. For the
inverse mapping $f^{-1}$ from relation (4) it follows:
$$f^{-1}: C_{x_1'}(\rho_0+2R))\rightarrow C_{x_1}(\rho_0+2R).  $$
Let us again connect points  $x_1'$ and  $y''$ with a $2R$-chain.
Moving along the chain by means of the same argument we get
$f^{-1}(y'')=x''\in X $. Finally, the mapping $f$ moves some point
$x''$ into the a priori chosen point  $y''$: $f(x'')=y''$.

Now we take a group $G_i=<f>$, where $f$ are all isometries of
$Iso(d)$ which superpose the clusters  $C_x(\rho_0)$ and
$C_{x'}(\rho_0)$,  $x, x'\in X_i$. We have proved that $G_i$
belongs to the group $G:=$Sym$(X)$ (i.e. $G_i\subseteq G$)  and
the group $G_i$ operates on $X_i$ transitively. In order to
complete the proof of Lemma 3.3 one needs to show that
$G_i\supseteq G$ for each $i\in[1,m]$. Indeed, we  show that if
$g\in G$, then $g\in G_i$. It is the case because $g$ is a
symmetry of $X$ and, hence, moves any point $y\in X_i$ and its
$\rho_0$-cluster  $C_y(\rho_0)$ into a point $g(y)\in X_i$ and the
cluster $C_{g(y)}(\rho_0)$, respectively. By the above-proved, the
symmetry $g$ belongs to the group $G_i$, i.e. $G_i\supseteq G$.
So, we proved that $G_i=$Sym$(X)$.
\endproof

\medskip
So, we have proved that a Delone set $X$ with conditions 1) and 2)
of Theorem 5 is partitioned into the union of  $m$ disjoint
discrete sets $X_i$ such that  each subset $X_i$ is a
$Sym(X)$-orbit of some point. In order to prove Theorem 5 we need
to make sure that $Sym(X)$ is a crystallographic group. A  proof
of  this fact is divided into two lemmas:

\begin{lemma}
Assume that   a group  $G\subset Iso(d)$ is such that for some point
$x\in \mathbb{R}^d$ its orbit $G\cdot x$ is a Delone set, then
$G$ is a crystallographic group.
\end{lemma}

\begin{lemma} In a partition $X=\bigsqcup_i^m X_i$  each  subset $X_i$ is a  Delone set.
\end{lemma}

From these two lemmas it  follows that $X$ is a crystallographic
orbit of finite set of $m$ points, i.e. that $X$ is a crystal with
respect to the group $G$. In fact, let us take a  Delone subset
$X_i$, which exists by lemma 3.5, and  a  group $G=G_i=$Sym$(X)$
generated by all possible isometries $f$ as in Lemma 3.3. Since
$G\cdot x=X_i$ where $x\in X_i$, by Lemma 2.4 the group $G$ is
crystallographic. Let a finite point set $X_0$ consist of $m$
representatives of subsets $X_i$, $i\in [1,m]$. Then $ X=G\cdot
X_0$,   We get that that $X$ is a crystallographic orbit of the
finite set $X_0$, i.e.  $X$ is a crystal. This completes the proof
of the local criterion.

Now we prove the last two lemmas.

\smallskip
\textbf{Proof} of Lemma 3.4. Let  $X:=G\cdot x$, (in  Lemma 3.4
$G$ is not assumed necessarily to be the full symmetry group
Sym$(X)$ of $X$, i.e. $G\subseteq Sym(X)$). Let   ${\cal V}$
denote the Voronoi tiling of space $\mathbb{R}^d $ with respect to
  $X$. A Voronoi domain  $V_x$ for the point  $x\in X$ is a cell
of the tiling $\cal V$.  $V_x$ is a   convex $d$-polytope with a
finite number of facets. This number of facets can be easily
estimated from above in terms of parameters  $r$ and $R$.
Therefore the symmetry group Sym$(V_x)$ of  $V_x$ is also finite.
Moreover, the order |Sym$(V_x)$| of this group can be also bounded
from above depending on the same parameters $r$ and $R$.

Every symmetry of the Delone set $X$ leaves the Voronoi  tiling
${\cal V}$ invariant. Therefore, since the group  $G$ operates on
the $X$ transitively, this group also operates transitively on the
set of all cells of the tiling ${\cal V}$. It is obvious that the
following inclusions are
 true:  $G\subseteq $Sym$(X)\subseteq$Sym$({\cal V}).$

The orbit Sym$({\cal V})\cdot y$ of any point $y\in \mathbb{R}^d$   is a discrete set
because the orbit Sym$({\cal V})\cdot y$ and Voronoi polytope $V_x$   intersect in  a finite set:
$$ |\hbox{Sym}({\cal V})\cdot y\cap V_x|\leq |\hbox{Sym}(V_x)|<c(r,R,d).$$
Therefore Sym$(X)$ is discrete.
Since $G$ is a  subgroup of Sym$(X)$, it is also discrete.

As for the fundamental domain $F(G)$, it can be chosen as $V_x/\hbox{stab}(x)$ where stab$(x)$ is the
stabilizer of the point $x$ in
$G$. In particular, if the stabilizer is trivial then the fundamental domain  $F(G)$ is the Voronoi polytope $V_x$.
 Thus,   the fundamental domain $F(G)$ is compact and hence  $G$ is a crystallographic group.\endproof

\textbf{Proof} of Lemma 3.5. First of all,  we note that since $X$
is a Delone set with parameter  $r$, any  subset $X_i$ fulfils the
(r)-condition of Delone set with a parameter $r'$, where $r'\geq
r$.

We suppose  that there is a  subset $X_i$ which does not satisfy
the second condition of Delone set for any finite value  $R'$. In
this case there is an infinite sequence of balls  $B_1,B_2,\ldots
, B_k,\ldots $ empty of points of $X_i$ with infinitely increasing
radii: $R_1<R_2<\cdots<R_k<\dots \rightarrow \infty$. Since the
set $X_i$ is discrete, each of these balls  $B_k$ can be moved so
that on its boundary there is a point  $x_k\in X_i$. Since all the
points  $x_k\in X_i$ belong to the $G$-orbit, one can move every
point  $x_k$ along with the ball $B_k$ to some given point $x\in
X_i$ by means of an appropriate isometry $f_k\in G$. Thus, one can
assume that the point $x$ is on
 the boundary of an empty ball
$B'_k$ of radius $R_k$ for every $k=1,2,\ldots $. Let
$\textbf{n}_k$ denote a unit vector
$$\textbf{n}_k:=\frac{1}{|xO_k|}\overrightarrow{xO_k},$$
where $O_k$ is the center of ball $B'_k$. Now we select from
sequence $\{\textbf{n}_k\}$ a subsequence
$\textbf{n}_{k_j}\rightarrow \textbf{n} $ that  converges to a
unit vector $\textbf{n}$.

Let  $\Pi$ be a hyperplane through point $x$ orthogonal to a
vector $ \textbf{n}$,  $\Pi^+$  the \textbf{open} half-space where
the normal vector $\textbf{n}$  points in. The half-space $\Pi^+$
contains no points of $X_1$. In fact, given a point  $z\in \Pi^+$,
in the subsequence of balls  with infinitely increasing radii one
can find a ball $B_{k_j}$ which contains  point $z$. Since all
balls $B_{k_j}$ are empty of points of $X_i$,
 $z\notin X_i$.

Thus, all points of  $X_i$ are in a  \textbf{closed }
half-space $\Pi^-$. We admit that all the points of $X_i$
can lay on the hyperplane $\Pi$ itself. Since $X$ is a Delone set, the open half-space  $\Pi^+$  also contains points of
$X$.

Given   $i\in[1,m]$ and $x\in X_i$, we take   $\j\neq i$ and
choose in $X_j$ a point $y$ closest to $x\in X_i$. One should note
that since $X_j$ is a discrete set, such a point $y$ does  exist.
Generally, there can be several but always finitely many closest
points. Let us denote  $\delta (x,X_j):=\min_{y\in Y_j}|xy|$.
Since $G$ operates transitively on both sets $X_i$ and $X_j$, the
minimum $\delta (x, X_j)$ does not depend on the choice of point
 $x\in Y_i$, i.e. for any other point $x'\in X_i$ there is a point  $y'\in X_j$ with condition
$\delta(x',X_j):=\delta_{ij}.$ It is clear also that
$\delta_{ij}=\delta_{ji}.$
Now  one can denote
$$\delta_i:=\max_{j\in [1,m], j\neq i}\delta_{ij}.$$
It is clear that  for every $j\in [1,m]$, $j\neq i$, and  $\forall
y\in X_j$ there is a point $x$ of $X_i$ at distance from $y$ of no
bigger  than $\delta_i$.

Therefore, since the subset $X_i$ is supposed to be located  in the closed half-space $\Pi^-$
 the whole  set $X$ is located in a half-space  $(\Pi+\delta_i
\textbf{n})^-$ determined by hyperplane  $\Pi+\delta_i \textbf{n}$.
The obtained contradiction
to the  $R$-condition of the  Delone set $X$ completes a proof of Lemma 3.5. \endproof

\section{Proof of Theorem 4}
At first we will prove Theorem 4. This theorem  easily implies
Theorems 1 and 2. The proof of Theorem 3 is based, in part on
Theorem 2. We note that in  Theorem 4 the locally antipodal Delone
sets $X$ and $Y$ a priori are not required  to be sets of finite
type.

Let us take an arbitrary point $x\in X$ and define the
\emph{distance spectrum} at the $x$ as the taken in ascending
order  set of distances between the point $x$ and the other points
of $X$:
$$\Re_{x}: =
\{\rho\in R_+| \exists x'\in X,\, |xx'|=\rho\}.$$ By the
$r$-condition for $X$ the spectrum $\Re_x$ is discrete and has no
limit points (with the exception of $\infty$) for any given $x\in
X$. Now we consider  the union $\cup_{x\in X}\Re_x$ over all $x\in
X$. It is easy to see that the union $\cup_{x\in X}\Re_{x}$ of the
spectra over all points of $X$  is a discrete set with no proper
limit point if and only if the Delone set $X$ is of finite type.

Recall  conditions for  Delone sets $X$ and $Y$: $x\in X\cap Y$
and $C_x(2R)=C'_x(2R)$. $C'_y(\rho)$ stands for the $\rho$-cluster
in the set $Y$). We  take the point $x$ and the two distance
spectra $\Re_x=\{\rho_1<\rho_2<\ldots \}$ and
$\Re'_x=\{\rho_1',\rho_2',\ldots \}$ in $X$ and in $Y$,
respectively, and prove the total coincidence of the spectra
$\Re_x$ and $\Re '_x$   and the sets $X$ and $Y$  by induction
along  numbers $k$  of  distance sequences $\rho_k$ and $\rho_k'$.

Due to the condition $C_x(2R) = C'_x(2R)$, some initial portions
of the spectra $\Re_x$ and $\Re'_x$ coincide. Assume that we have
already proved the equality of the first $k$ distances
$\rho_1=\rho'_1, \ldots , \rho_k=\rho'_k$ in the spectra and the
coincidence of the clusters $C_x(\rho_k)=C'_x(\rho_k)$.

Now we prove that $\rho_{k+1}=\rho'_{k+1}$ and
$C_x(\rho_{k+1})=C'_x(\rho'_{k+1})$. Let $\rho_{k+1}\leq
\rho'_{k+1}$, then the ball $B_x(\rho_{k+1})$ has on its boundary
at least one  point $x_1\in X$, $|xx_1|=\rho_{k+1}$ (see Fig. 4).
Let $z\in \mathbb{R}^d$ be such that $z\in [xx_1]$ and the length
$|zx_1|=R$. We note that point $z$, generally,   does not belong
to $X$. The ball $B_z(R)$ centered at $z$  of radius $R$
 touches the sphere $\partial B_x(\rho_{k+1})$ at point $x_1$. Now let $h$ be the  homothety  with the center $x_1$ and the
coefficient 2 and $B_{z'}(2R):=h(B_z(R))$ (see Fig. 4). It is
obvious that we have
$$B_z(R)\subset B_{z'}(2R)\subset B^o_x(\rho_{k+1})\cup \{x_1\}. \eqno(2)
 $$
Here $B^o$ means an open ball.

\begin{center}
\begin{figure}[!ht]
\hbox{\includegraphics[width=8cm]{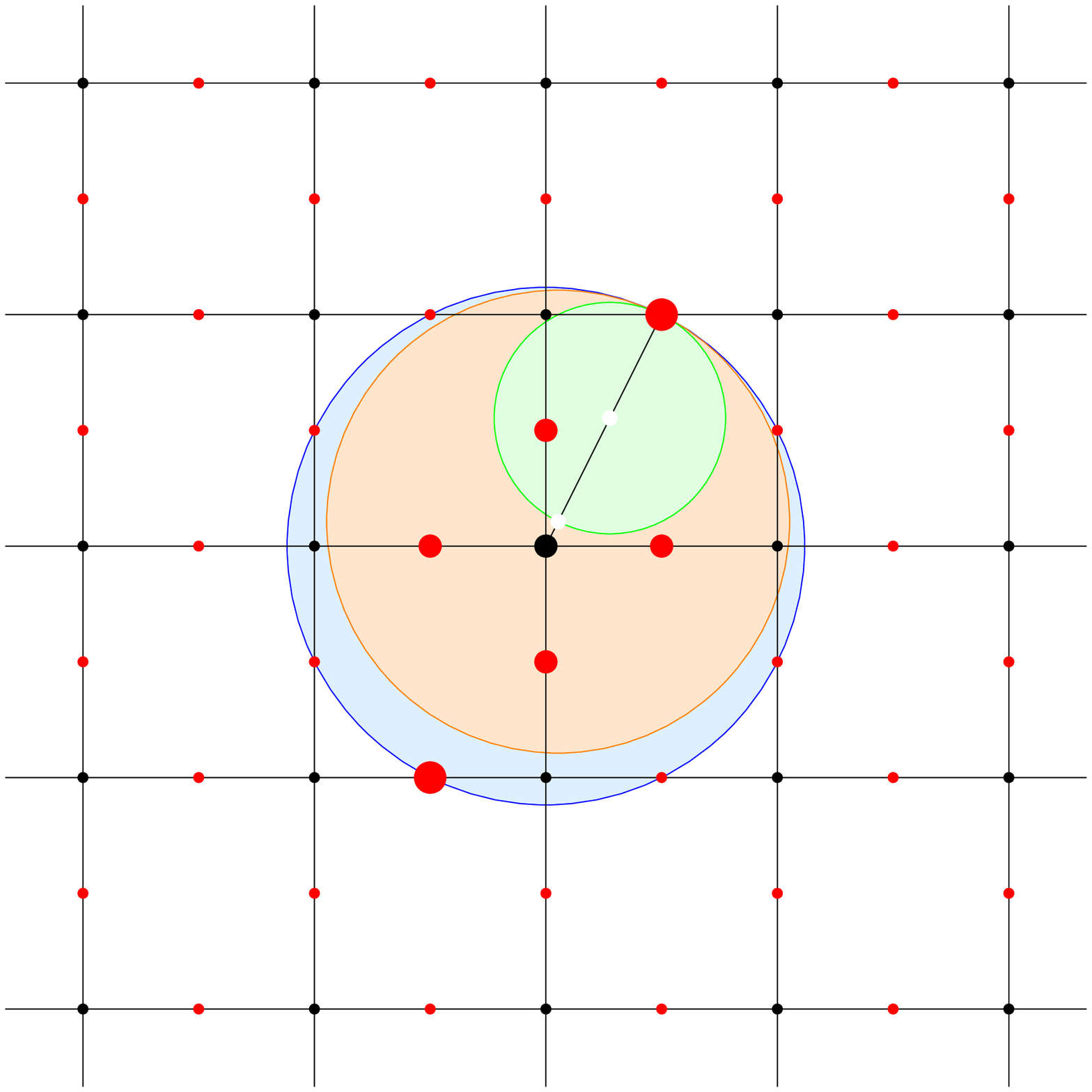}
\raise5.0cm\llap{$z$\kern3.1cm}%
 \raise4.4cm\llap{$z'$\kern3.3cm}%
 \raise4.0cm\llap{$x$\kern3.9cm}%
\raise6.0cm\llap{$x_1$\kern3.0cm}%
\raise5.0cm\llap{$x_2$\kern3.9cm}%
\raise4.1cm\llap{$x_3$\kern4.6cm}%
\raise0.0cm\llap{Fig. 4\kern3.4cm}%
}
\end{figure}
\end{center}

By the $R$-condition, in $B_z(R)$ there is at least one  point
$x_2\in X$, $x_2\neq x_1$. Since $x_1$ is the only point of the
ball $B_z(R)$ which is located  on the boundary $\partial
B_x(\rho_{k+1})$, all other points of $B_z(R)$, including the
point $x_2$, lay in the interior of $B_x(\rho_{k+1})$. Therefore
$|xx_2|<|xx_1|=\rho_{k+1}$, i.e. $|xx_2|\leq \rho_k|$, and,  by
the induction assumption,    $x_2\in X\cap Y$. Since $|x_1x_2|\leq
2R$ the point $x_1$ belongs to the cluster $ C_{x_2}(2R)$.

 Since the cluster $C_{x_2}(2R)$ is  antipodal  about $x_2$,
 in this cluster there is a point $x_3$ which is
 antipodal to $x_1$.
 We recall that the coefficient of the homothety $h$ equals 2, hence
$$x_3\in B_{z'}(2R)\subset B^o_x(\rho_{k+1})\cup\{x_1\}.$$
 Therefore, $|xx_3|\leq \rho_k$, as well as $|xx_2|\leq \rho_k$,
 and,  by the inductive assumption, $x_2,x_3\in Y$.
Since $|x_2x_3|\leq 2R$ ,
 we have that
$x_3\in C'_{x_2}(2R)$. Now, since $x_1$ is antipodal to $x_3$
about $x_2$ and the cluster $C'_{x_2}(2R)$ is antipodal about
$x_2$, $x_1$ also belongs to $C'_{x_2}(2R)$. Hence we have also
$x_1\in Y$. This inclusion is true for any $x_1'\in X$ with
$|xx_1'|=\rho_{k+1}$. Thus, it has been  proved that if
$\rho_{k+1}\leq \rho'_{k+1}$ we actually have
$\rho_{k+1}=\rho'_{k+1}$. Therefore we just  proved
$C_x(\rho_{k+1})\subseteq C'_x(\rho_{k+1})$. However, in the case
$\rho_{k+1}=\rho'_{k+1}$ one can also take any point $y_1\in Y$
with $|xy_1|=\rho_{k+1}$ and by the  same argument
 prove that $y_1\in X$. Thus, the inductive step is established: one has proved that $C_x(\rho_{k+1})= C'_x(\rho_{k+1})$.
  \endproof


\section{Proofs of Theorems 1,2, and  3}

Theorems 1 and 2 easily follow from Theorem 4.

\medskip
\textbf{Proof of Theorem 1}.  By the requirement $N(2R)=1$, for
any $x'$ and $x\in X$ there is an isometry $g$ such that $g(x')=x$
and
 $g(C_{x'}(2R))=C_x(2R)$.  Let us denote $Y:=g(X)$. We have two local antipodal sets $X$ and $Y$ such that
$X\cap Y \supseteq C_x(2R)$. By Theorem 4, the relationship
$C_x(2R)=C'_x(2R)$ implies $X=Y$, i.e. $g$ is a symmetry of $X$.
Thus Sym$(X)$ possesses  transitive symmetry group.\endproof

\medskip
\textbf{Proof of Theorem 2}. $X$ is a locally antipodal Delone
set. Let $\sigma_x$ be the  inversion  about a point $x$ such that
$\sigma_x: C_x(2R)\rightarrow C_x(2R)$. Let us denote $Y:=\sigma
(X)$. Then we have again two sets $X$ and $Y$ with a $2R$-cluster
$C_x(2R)$ in common. By Theorem 4 we have $X=Y$, i.e. the
inversion $\sigma_x$ leaves the set $X$ invariant.
\endproof

\medskip
\textbf{Proof of Theorem 3}

Given a locally antipodal set $X\subset \mathbb{R}^d $, let
$\Lambda$ be a set of vectors $\lambda$ such that $X+\Lambda =X$.
Since $X$ is discrete the vector set $\Lambda$ is a lattice.

Now we show that the lattice $\Lambda$ is a lattice of  rank $d$.
Indeed, let $\sigma_x$ and $\sigma_{x'}$ be inversions of clusters
$C_x(2R)$ and $C_{x'}(2R)$ at points $x$ and $x'$, respectively.
Then,  by Theorem 2, they are both symmetries of the whole $X$. On
the other hand, the superposition $\sigma_x\circ\sigma_{x'}$ is a
translation by the vector $2(x-x')$. Since the set $X$ is a Delone
set, the  translational group
 $\Lambda$ generated by all possible $2(x-x')$, where $x,x'\in X$, is a lattice of the rank $d$.
$\Lambda $ is the maximum lattice to leave $X$ invariant, i.e.
such that $X+\Lambda=X$.

The Delone set  $X$  is the union of finitely many lattices which
are congruent and parallel to the lattice $\Lambda$, i.e.
$$X = \bigcup\limits_{i = 1}^n (x_i + \Lambda).$$
As proved above, for  $i = 2, 3, \ldots, n$ we have  $x_i - x_1
\in \Lambda / 2$.

By putting  $x := x_1$, $\lambda_i/2 = x_i - x_1$ ($i =
1, 2, \ldots, n$) we come to:
$$X=\bigsqcup\limits_{i=1}^n(x+\lambda_i/2+\Lambda), \hbox{ where } \lambda_i\in \Lambda. \eqno(6)$$
Now, if $\lambda_i\equiv \lambda_j\mod 2$, i.e. if
$\lambda_i-\lambda_j=2\lambda,$ where $\lambda\in \Lambda$, then
subsets $x+\lambda_i + \Lambda$ and $x+ \lambda_j + \Lambda$
obviously coincide. Therefore in (6) $n\leq 2^d$. Moreover, the
value $n$ cannot be equal to $2^d$ because in this case
$X=\Lambda/2$ and hence $X+\Lambda/2=X$. This contradicts the
assumption  that $\Lambda$ is the  maximum  lattice with the
condition $X+\Lambda$. So, $n\leq 2^d-1$. \endproof

\vskip 0.3cm The author is very grateful to Nikolay Andreev
(Moscow) for making drawings and to Andrey Ordine (Toronto) for
his help in editing the English text.

\section*{References}.

\medskip

\noindent [1] B. Delaunay, Sur la sphere vide. A la memoire de
Georges Voronoi.
Bulletin de l'Academie des Sciences de l'URSS. Classe des sciences mathematiques et na, 1934, Issue 6,  Pages 793–800.

\noindent [2] B.N. Delone, Geometry of positive quadratic forms, Uspekhi Matem. Nauk, 1937, 3,  16–62 (in Russian).

\noindent [3] A. Sch\"oflies, Kristallsysteme und
Kristallstruktur, Leipzig,  1891 - Druck und verlag von BG Teubner

\noindent [4] L.Bieberbach, \"Ueber die Bewegungsgruppen des
n-dimensionalen Euklidischen R\"aumes I, Math. Ann. 70 (1911),
207-336; II, Math. Ann. 72 (1912), 400-412.

\noindent [5] B. N. Delone, N.P. Dolbilin, M.I. Stogrin, R.V. Galiuilin,
A local criterion for regularity of a system of points, Soviet Math. Dokl., 17, 1976, 319–322.

\noindent [6] N. P. Dolbilin, M. I. Shtogrin, A local criterion
for a crystal structure, Abstracts of the IXth All-Union
Geometrical Conference, Kishinev, 1988, p. 99 (in Russian).

\noindent [7] N.P.~Dolbilin, A Criterion for   crystal and locally
antipodal Delone sets. Vestnik Chelyabinskogo Gos. Universiteta,
2015, 3 (358), 6--17 (in Russian).

\noindent [8] N.P.~Dolbilin, A.N.~Magazinov, Locally antipodal
Delauney Sets, Russian Math. Surveys.70:5 (2015), 958-960.

\noindent [9] E.S.~Fedorov, Elements of the Study of Figures, Zap.
Mineral. Imper. S.Peterburgskogo Obschestva, 21(2), 1985, 1-279.

 \noindent [10] R.~Feynman, R.~Leighton, M.~Sands, Feynman Lectures on Physics, Vol. II, Addison-Wesley, 1964.

\noindent [11] N.P.~Dolbilin, J.C.~Lagarias, M.~Senechal,
Multiregular point systems. Discr. and Comput. Geometry, 20, 1998,
477–498.

\noindent [12]     N.~Dolbilin, E.~Schulte, The local theorem for
monotypic tilings, Electron. J. Combin., 11:2 (2004), Research
Paper 7 , 19 pp.

\noindent [13]  N.~Dolbilin, E.~Schulte, A local characterization
of combinatorial multihedrality in tilings, Contrib. Discrete
Math., 4:1 (2009), 1–11.

\noindent [14] N.~Dolbilin, Regular systems in 3D space. Chebyshev
Sbornik, (2016) (in print).

\noindent [15] D.~Schattschneider, N.~Dolbilin, One corona is
enough for the Euclidean plane. Quasicrystals and Discrete
Geometry, Fields Inst. Monogr., 10, American Math. Soc. Providence
RI, 1998, 207-246.

\end{document}